\documentclass[11pt]{article}
\usepackage{amsfonts}
\newcommand{\be}{\begin{equation}}
\newcommand{\ee}{\end{equation}}
\newcommand{\bea}{\begin{eqnarray}}
\newcommand{\eea}{\end{eqnarray}}
\newcommand{\beaa}{\begin{eqnarray*}}
\newcommand{\eeaa}{\end{eqnarray*}}
\newcommand{\eq}[1]{(\ref{#1})}

\newcommand{\eps}{\varepsilon}
\newcommand{\vphi}{\varphi}
\newcommand{\tea}{\vartheta}

\newcommand{\Ricc}{\mathbf{Ricci}}
\newcommand{\metric}{\mathbf{g}}
\newcommand{\energy}{{\mathbf{T}}}
\newcommand{\zero}[1]{\stackrel{0}{#1}}
\newcommand{\one}[1]{\stackrel{1}{#1}}
\newcommand{\two}[1]{\stackrel{2}{#1}}
\newcommand{\mone}[1]{\stackrel{-1}{#1}}
\newcommand{\mtwo}[1]{\stackrel{-2}{#1}}
\begin{document}

\renewcommand{\theequation}{\thesection.\arabic{equation}}

\title{Large Amplitude Gravitational Waves}
\author{G.~Al\`{\i} \\ Institute for Applications of Mathematics \\
Consiglio Nazionale delle Ricerche, Napoli \and J.~K.~Hunter
\\ Department of Mathematics and \\ Institute of Theoretical Dynamics \\
University of California at Davis}

\date{November 25, 1998}

\maketitle

\begin{abstract}

We derive an asymptotic solution of the
Einstein field equations which
describes the propagation of a thin, large
amplitude gravitational wave into a curved space-time.
The resulting equations have the same form
as the colliding plane wave equations without
one of the usual constraint equations.
\end{abstract} 

\section{Introduction}
\setcounter{equation}{0}

Gravitational waves are one of the most important
features of Einstein's general theory of relativity.
The Einstein field equations
are highly nonlinear, and a question of fundamental
interest is how nonlinearity affects the propagation
of gravitational waves. 

Small amplitude gravitational waves are well described by the
linearized Einstein equations which completely
neglect nonlinear effects. Large amplitude unidirectional gravitational
plane waves are described by the exact
Brinkmann-Rosen solution of the vacuum Einstein equations \cite{Br,LL}.
Despite the nonlinearity of the Einstein equations,
a gravitational plane wave
propagates into flat space-time without distortion, and
there are no dynamic nonlinear effects.

The simplest situation in which
nonlinear effects are significant
is when a large amplitude
gravitational wave propagates into a curved space-time.
An important special case is when the space-time ahead
of the wave is that of a counter-propagating gravitational
plane wave. The resulting space-time has
a two-parameter family of spacelike isometries, and the
metric is given by the exact colliding plane wave solution
of the vacuum Einstein equations \cite{KP,Sz1,Sz2,Gr}.
In the case of more general space-times ahead of the wave,
exact solutions do not exist.

In this paper, we derive an asymptotic solution of the Einstein equations
which describes the propagation of a thin, large-amplitude,
pulse-like gravitational wave into a general curved space-time.
The solution applies provided that
the metric varies much more rapidly inside
the wave than on either side of the wave.
As a result, the wave can be approximated
locally by a nonlinear plane wave which
is distorted as it propagates into the
curved space-time. For plane-polarized waves,
the asymptotic solution is given by equations \eq{asysol}, \eq{plane-polar},
\eq{thetaconstraint}, and \eq{apc1}--\eq{apc3} below.
For non-polarized waves, the asymptotic solution
is given by \eq{asysol}, \eq{Brp}, \eq{np0},
and \eq{np1}--\eq{np4}. The asymptotic equations
consist of the colliding plane wave equations
without one of the usual constraint equations.
The colliding plane wave equations are
therefore canonical equations for
nonlinear gravitational waves, and they describe
a much larger class of solutions than the ones with exact
plane-wave symmetry.

The nonlinearity of the asymptotic equations may
lead to the development of a space-time singularity.
A plane gravitational wave propagating into flat space-time
does not steepen. Consequently,
the mechanism of singularity formation in gravitational
waves differs from the nonlinear steepening of waves in quasilinear
hyperbolic systems which leads to the formation of shocks.
Instead, the singularity formation is
caused by the mutual focusing of the gravitational
wave and the curved space-time into which it propagates.
A second nonlinear effect described by the asymptotic equations
is the permanent distortion of space-time by the passage
of a gravitational wave. A curved pulse
generates a backscattered gravitational wave which propagates into
the space-time behind it.

In Section 2, we briefly summarize the exact colliding plane
wave solution of the Einstein equations. In Section 3, we give
an overview of the asymptotic expansion. In Section 4, we write
out expansions of the metric components,
the connection coefficients, and the Ricci
curvature components. In Section 5, we construct a coordinate system in
which the metric adopts its simplest form. In Section 6, we complete
the derivation of the asymptotic equations. In Section 7, we show that
the same equations follow from an expansion of the variational principle
for the Einstein equations. In Section 8, we explain how to derive
boundary conditions for the asymptotic equations, and in Section
9 we consider some specific physical examples.

\section{Colliding plane waves}
\setcounter{equation}{0}

The vacuum Einstein field equations imply that
\be
\Ricc = 0,\label{einsteineq}
\ee
where $\Ricc$ is the Ricci tensor associated with
the metric tensor $\metric$.
The plane-plolarized, colliding plane wave solution of 
\eq{einsteineq} is given by
\be
\metric=-2e^{-M}du\,dv+e^{-U}(e^{V}dy^2+e^{-V}dz^2),\label{Br}
\ee
where the functions $M(u,v)$, $U(u,v)$, $V(u,v)$ satisfy the
colliding plane wave equations,
\bea
U_{uv} &=& U_{u}U_{v},
\label{pc1}
\\
V_{uv} &=& \frac{1}{2}\left(U_{u}V_{v}+U_{v}V_{u}\right),
\label{pc2}
\\
M_{uv} &=& \frac{1}{2}\left(-U_{u}U_{v}+V_{u}V_{v}\right),
\label{pc3}
\\
U_{uu} &=& \frac{1}{2}\left(U_{u}^2+V_{u}^2\right)-U_{u}M_{u},
\label{pc4}
\\
U_{vv} &=& \frac{1}{2}\left(U_{v}^2+V_{v}^2\right)-U_{v}M_{v}.
\label{pc5}
\eea
Equations (\ref{pc1})--(\ref{pc3}) are wave equations
for $M$, $U$ and $V$ in characteristic coordinates $(u,v)$.
Equations (\ref{pc4})--(\ref{pc5})
are constraints which are preserved by (\ref{pc1})--(\ref{pc3}).
To specify a unique solution,
the wave equations can be supplemented
by characteristic initial data for $M$, $U$, $V$
on the lines $u=0$ and $v=0$ which satisfy the
appropriate constraint equations.

The metric which describes the collision of
non-polarized plane waves is
\[
\metric = -2e^{-M}du\,dv
+e^{-U}(e^{V}\cosh W dy^2-2\sinh W dy\,dz+e^{-V}\cosh W dz^2),
\]
where the functions $M(u,v)$, $U(u,v)$, $V(u,v)$, $W(u,v)$ satisfy
\bea
U_{uv} &=& U_{u}U_{v},
\label{npc1}
\\
V_{uv} &=& \frac{1}{2}(U_{u}V_{v}+U_{v}V_{u})
-(V_{u}W_{v}+V_{v}W_{u})\tanh W,
\label{npc2}
\\
W_{uv} &=& \frac{1}{2}(U_{u}W_{v}+U_{v}W_{u})
+V_{u}V_{v}\sinh W\cosh W,
\label{npc3}
\\
M_{uv} &=& \frac{1}{2}(-U_{u}U_{v}+V_{u}V_{v}\cosh^2 W+W_{u}W_{v}),
\label{npc4}
\\
U_{uu} &=& \frac{1}{2}(U_{u}^2+V_{u}^2\cosh^2 W+W_{u}^2)-U_{u}M_{u},
\label{npc5}
\\
U_{vv} &=& \frac{1}{2}(U_{v}^2+V_{v}^2\cosh^2 W+W_{v}^2)-U_{v}M_{v}.
\label{npc6}
\eea
When $W=0$, this solution reduces to the plane-polarized solution.
When all functions are independent of $v$, the solution reduces to the
Rosen form of the exact unidirectional plane wave solution.

\section{Overview of the expansion}
\label{plane}
\setcounter{equation}{0}

In this section, we outline the main ideas of the derivation
of the asymptotic solution.
For simplicity, we describe the case of plane-polarized waves.
The algebraic details are given in the following sections.

We consider metrics of the form
\bea
\metric &=& \metric\left(\frac{u(x)}{\eps}, x ; \eps\right),
\label{asysol}
\\
\metric\left(\theta,x;\eps\right) &=& \zero{\metric}\left(\theta, x\right)
+ \eps\one{\metric}\left(\theta, x\right) + O(\eps^2),
\nonumber
\eea
where $\eps$ is a small parameter and $u$ is a scalar-valued phase
function with $du \ne 0$. This ansatz corresponds to a
metric that varies rapidly and strongly in the $u$-direction.
The phase $u$ is a null function of the metric,
at least up to the order $\eps$. That is, it satisfies
\be
\metric^\sharp(du,du) = O(\eps^2),
\label{null}
\ee
where $\metric^\sharp$ is the contravariant form of the metric
tensor. The component form of this equation is written
out in \eq{nullco} below.

The scaled variable 
\be
\theta = \frac{u}{\eps}\label{deftheta}
\ee
is a ``stretched'' coordinate inside  the wave.
We assume that the derivatives of $\metric(\theta,x;\eps)$
with respect to $\theta$ decay to zero sufficiently quickly as
$\theta \to \infty$. Thus, the solution \eq{asysol} represents
a thin, pulse-like gravitational wave located
near the null surface $u = 0$. For example, if the metric is
independent of $\theta$ when $|\theta|$ is sufficiently large,
then the solution represents a thin ``sandwich'' wave
which separates slowly varying metrics on either side.

The Ricci tensor associated with the metric \eq{asysol}
has an expansion of the form
\be
\Ricc = \frac{1}{\eps^2}\mtwo{\Ricc}+ \frac{1}{\eps}\mone{\Ricc}
+ O(1).
\ee
At leading order in $\eps$, the Einstein equations \eq{einsteineq}
imply that
\[
\mtwo{\Ricc} = 0.
\]
This equation is a nonlinear, second order ordinary
differential equation in $\partial_\theta$ for
the leading order term of the metric in which
the ``slow'' variables $x$ occur as parameters. We write it
symbolically as
\be
N(\partial_\theta^2)\left[\zero{\metric}\right] = 0.\label{equzero}
\ee
In suitable coordinates $(u,v,y,z)$, a solution of this
equation is the plane-polarized plane wave metric
\be
\zero{\metric}= -2e^{-M}du\, dv
+e^{-U}\left(e^{V} dy^2 + e^{-V}dz^2\right),
\label{plane-polar}
\ee
where  $M$, $U$, $V$ are functions of $(\theta,v,y,z)$.
For a metric of the form \eq{plane-polar}, equation \eq{equzero}
reduces to the $\theta$-constraint equation,
\be
U_{\theta\theta} = \frac12\left(U_\theta^2 + V_\theta^2\right)
-U_{\theta}M_{\theta}.\label{thetaconstraint}
\ee

At the next order in $\eps$, the Einstein equations imply that
\[
\mone{\Ricc} = 0.
\]
This is a linear equation for $\one{\metric}$ of the form
\be
L(\partial_\theta^2)\left[\one{\metric}\right] =
F(\partial_\theta, \partial_v,\partial_y,\partial_z)
\left[\zero{\metric}\right],\label{equone}
\ee
where $L$ is a second order linear ordinary differential
operator in $\partial_\theta$ acting on $\one{\metric}$,
with coefficients depending
on $\zero{\metric}$, and $F$ is a nonlinear partial differential
operator acting on $\zero{\metric}$.
The equations in \eq{equone} are not independent.
The requirement that \eq{equone} can be solved for $\one{\metric}$
implies that $M$, $U$, and $V$ satisfy the equations
\bea
&&U_{\theta v} =
U_{\theta}U_{v},
\label{apc1}
\\
&&V_{\theta v} =
\frac{1}{2}\left(U_{\theta}V_{v}
+U_{v}V_{\theta}\right),
\label{apc2}
\\
&&M_{\theta v} =  
\frac12 \left(-U_\theta U_{v}
+ V_\theta V_{v}\right).
\label{apc3}
\eea
Equations \eq{apc1}--\eq{apc3} are
identical to the evolution equations \eq{npc1}--\eq{npc3}
for the exact colliding plane wave solution, with
$\theta=u/\eps$.  The leading order solution
satisfies the constraint equation \eq{thetaconstraint}
in the ``fast'' phase variable $\theta$, but it need not
satisfy the constraint equation \eq{pc5}
in the ``slow'' variable $v$. If the $v$-constraint
equation does not hold, then the asymptotic expansion
of the metric contains higher order terms which are
absent in the exact colliding plane wave solution.

Equation \eq{plane-polar} implies that
$\partial_v = - e^M \metric^\sharp\cdot du$.
Thus, $\partial_v$ is a vector on the
light cone which is tangent to the null surface $u=0$,
and the ``slow'' derivative with respect to $v$ which
appears in \eq{apc1}--\eq{apc3} is a derivative along
the bicharacteristic null geodesics associated with $u$.
The transverse variables $y$ and $z$ occur as parameters.
Therefore, in the short-wave limit considered here, the $(1+3)$-dimensional
field equations reduce to $(1+1)$-dimensional asymptotic
equations along the set of null geodesics associated with
the phase $u$. The parametric dependence of the solution on $y$
and $z$ allows the pulse to be compactly supported in the 
transverse directions, so that the wave need not have infinite
extent. Moreover, the asymptotic solution need not have
any special exact symmetries.

The asymptotic equations for non-polarized gravitational waves are
obtained in a similar way. They consist of the general colliding
plane wave equations \eq{npc1}--\eq{npc5} with
$u$ replaced by $\theta$. The $v$-constraint equation \eq{npc6}
is not required to hold.

Since the asymptotic equations follow from the
order $\eps^{-2}$ and order $\eps^{-1}$ components of the
field equations, the asymptotic solution remains valid
in the presence of matter with a slowly varying, order
one energy-momentum tensor, $\energy = \energy(x)$.

One subtle point in carrying out the expansion concerns the
choice of the phase function $u$. In order for \eq{equzero}
to have a nontrivial solution, the phase $u$ must be a null function
of the leading order metric, but $u$ need not be a null
function of the entire metric. However,
it follows from the analysis in Section \ref{gauge}
that we can use a transformation of the form
\be
u \to  \eps \Psi\left(\frac{u}{\eps},x;\eps\right)
\label{newu}
\ee
to choose a phase which satisfies \eq{null}. The asymptotic solutions
obtained with the use of the old and the new phases can be shown to be
equivalent. When the phase satisfies \eq{null}, variations
in the metric propagate along the null geodesics associated with the
phase, and the asymptotic equations adopt their simplest form.

\section{Expansion of the metric and the curvature}
\label{asymptotic}
\setcounter{equation}{0}

In this section, we write out expansions
of the metric components, the connection coefficients, and the
Ricci curvature components.

We use local coordinates $x^\alpha$ in which
\be
\metric=g_{\alpha\beta}dx^\alpha\,dx^\beta.
\label{metriccomp}
\ee
Here and below, Greek indices $\alpha,\beta,\mu,\nu,\dots$ take
on the values $0,1,2,3$.
We look for an expansion of the metric components as $\eps \to 0$
of the form
\bea
g_{\alpha\beta} &=& g_{\alpha\beta}\left(\frac{u(x)}{\eps},x;\eps\right),
\label{exp1}
\\
g_{\alpha\beta}\left(\theta,x;\eps\right) 
&=& \zero{g}_{\alpha\beta}\left(\theta,x\right)
+\eps \one{g}_{\alpha\beta}\left(\theta,x\right)
+O(\eps^2).
\nonumber
\eea
The contravariant metric components $g^{\alpha\beta}$ satisfy
\[
g^{\alpha\mu}g_{\mu\beta}=\delta^\alpha_\beta.
\]
Expansion of this equation in a power series in $\eps$
gives
\be
g^{\alpha\beta}
=\zero{g}{}\!^{\alpha\beta}
-\eps \one{g}{}\!^{\alpha\beta} +O(\eps^2).
\label{exp2}
\ee
In \eq{exp2}, $\zero{g}{}\!^{\alpha\beta}$ is the inverse
of $\zero{g}_{\alpha\beta}$, and
we use the leading order metric
components to raise indices, so that
\be
\one{g}{}\!^{\alpha\beta}=
\zero{g}{}\!^{\alpha\mu}\zero{g}{}\!^{\beta\nu}
\one{g}_{\mu\nu}.
\label{g1}
\ee
With this notation, the order $\eps$ term in the expansion
of the contravariant
metric component $g^{\alpha\beta}$ is $-\one{g}{}\!^{\alpha\beta}$, not
$\one{g}{}\!^{\alpha\beta}$.

In terms of the metric components, we have
\bea
\metric^\sharp(du,du) &=&  g^{\alpha\beta}
\frac{\partial u}{\partial x^\alpha}\frac{\partial u}{\partial x^\beta}
\nonumber
\\
&=&
\zero{g}{}\!^{\alpha\beta}\frac{\partial u}{\partial x^\alpha}
\frac{\partial u}{\partial x^\beta}
- \eps \one{g}{}\!^{\alpha\beta}\frac{\partial u}{\partial x^\alpha}
\frac{\partial u}{\partial x^\beta} + O(\eps^2).
\label{nullco}
\eea
Thus, the null condition \eq{null} holds provided that
\be
\zero{g}{}\!^{\alpha\beta}\frac{\partial u}{\partial x^\alpha}
\frac{\partial u}{\partial x^\beta} = 0, \qquad
\one{g}{}\!^{\alpha\beta}\frac{\partial u}{\partial x^\alpha}
\frac{\partial u}{\partial x^\beta} = 0. \label{null1}
\ee
The first condition  in \eq{null1} states that $u$ is a null function
of $\zero\metric$. The second condition
implies that the phase is a null function
of the perturbed metric, at least up to the first order in $\eps$.

The covariant components $R_{\alpha\beta}$
of the Ricci tensor are given by
\be
R_{\alpha\beta}= 
\frac{\partial\Gamma^\lambda{}_{\alpha\beta}}{\partial x^\lambda}
-\frac{\partial\Gamma^\lambda{}_{\beta\lambda}}{\partial x^\alpha}
+\Gamma^\lambda{}_{\alpha\beta}\Gamma^\mu{}_{\lambda\mu}
-\Gamma^\mu{}_{\alpha\lambda}\Gamma^\lambda{}_{\beta\mu},
\label{ricci}
\ee
where $\Gamma^\lambda{}_{\alpha\beta}$ are the connection coefficients
\be
\Gamma^\lambda{}_{\alpha\beta}=
\frac{1}{2}g^{\lambda\mu}
\left(
\frac{\partial g_{\beta\mu}}{\partial x^\alpha}
+\frac{\partial g_{\alpha\mu}}{\partial x^\beta}
-\frac{\partial g_{\alpha\beta}}{\partial x^\mu}
\right).
\label{connection}
\ee
>From \eq{deftheta}, the derivative of a function $f_{\alpha\beta}(\theta,x)$,
with respect to $x^\mu$ is given by
\be
\frac{\partial f_{\alpha\beta}}{\partial x^{\mu}}
 =
\frac{1}{\eps} f_{\alpha\beta,\theta}u_\mu +f_{\alpha\beta,\mu},
\label{derivexp}
\ee
where
\[
u_\mu = \frac{\partial u}{\partial x^{\mu}}, \quad f_{\alpha\beta,\theta}
=\left.\frac{\partial f_{\alpha\beta}}{\partial\theta}\right|_{x},
\quad
f_{\alpha\beta,\mu}
=\left.\frac{\partial f_{\alpha\beta}}{\partial x^\mu}\right|_{\theta}.
\]
We use \eq{exp1}, \eq{exp2}, and \eq{derivexp} in \eq{ricci}
and \eq{connection} and expand the result with respect to $\eps$.
After some algebra, we find that
\bea
&&\Gamma^\lambda{}_{\alpha\beta}=
\frac{1}{\eps}\mone{\Gamma}{}\!^\lambda{}_{\alpha\beta}
+\zero{\Gamma}{}\!^\lambda{}_{\alpha\beta}
+O(\eps),
\nonumber
\\
&&R_{\alpha\beta}=
\frac{1}{\eps^2}\mtwo{R}_{\alpha\beta}
+\frac{1}{\eps}\mone{R}_{\alpha\beta}
+O(1),
\label{Rexp}
\eea
where
\bea
\mone{\Gamma}{}\!^\lambda{}_{\alpha\beta}
&=&
\frac{1}{2}\zero{g}{}\!^{\lambda\mu}
\left(
\zero{g}_{\beta\mu,\theta}u_\alpha
+\zero{g}_{\alpha\mu,\theta}u_\beta
-\zero{g}_{\alpha\beta,\theta}u_\mu
\right),
\nonumber
\\
\zero{\Gamma}{}\!^\lambda{}_{\alpha\beta}
&=&
\frac{1}{2}\zero{g}{}\!^{\lambda\mu}
\left(
\zero{g}_{\beta\mu,\alpha}
+\zero{g}_{\alpha\mu,\beta}
-\zero{g}_{\alpha\beta,\mu}
\right)
\nonumber
\\
&&
+\frac{1}{2}\zero{g}{}\!^{\lambda\mu}
\left(
\one{g}_{\beta\mu,\theta}u_\alpha
+\one{g}_{\alpha\mu,\theta}u_\beta
-\one{g}_{\alpha\beta,\theta}u_\mu
\right)
\nonumber
\\
&&
-\frac{1}{2}\one{g}{}\!^{\lambda\mu}
\left(
\zero{g}_{\beta\mu,\theta}u_\alpha
+\zero{g}_{\alpha\mu,\theta}u_\beta
-\zero{g}_{\alpha\beta,\theta}u_\mu
\right),
\label{connricciexp}
\\
\mtwo{R}_{\alpha\beta}
&=&
\mone{\Gamma}{}\!^\mu{}_{\alpha\beta,\theta}u_\mu
-\mone{\Gamma}{}\!^\mu{}_{\beta\mu,\theta}u_\alpha
+\mone{\Gamma}{}\!^\mu{}_{\alpha\beta}
\mone{\Gamma}{}\!^\nu{}_{\mu\nu}
-\mone{\Gamma}{}\!^\mu{}_{\alpha\nu}
\mone{\Gamma}{}\!^\nu{}_{\beta\mu},
\nonumber
\\
\mone{R}_{\alpha\beta}
&=&
\mone{\Gamma}{}\!^\mu{}_{\alpha\beta,\mu}
-\mone{\Gamma}{}\!^\mu{}_{\beta\mu,\alpha}
+\zero{\Gamma}{}\!^\mu{}_{\alpha\beta,\theta}u_\mu
-\zero{\Gamma}{}\!^\mu{}_{\beta\mu,\theta}u_\alpha
\nonumber
\\
&&
+\mone{\Gamma}{}\!^\mu{}_{\alpha\beta}
\zero{\Gamma}{}\!^\nu{}_{\mu\nu}
+\zero{\Gamma}{}\!^\mu{}_{\alpha\beta}
\mone{\Gamma}{}\!^\nu{}_{\mu\nu}
-\mone{\Gamma}{}\!^\mu{}_{\alpha\nu}
\zero{\Gamma}{}\!^\nu{}_{\beta\mu}
-\zero{\Gamma}{}\!^\mu{}_{\alpha\nu}
\mone{\Gamma}{}\!^\nu{}_{\beta\mu}.
\nonumber
\eea

The component form of the field equations \eq{einsteineq} is
\be
R_{\alpha\beta} = 0.\label{einsteineqco}
\ee
Using \eq{Rexp} in \eq{einsteineqco} and equating coefficients
of $\eps^{-2}$ and $\eps^{-1}$ to zero, we get that
\bea
\mtwo{R}_{\alpha\beta}&=&0,
\label{pert1}
\\
\mone{R}_{\alpha\beta} &=& 0.
\label{pert2}
\eea
In order to solve these equations, we first use a coordinate
transformation to simplify the form of the metric.

\section{Coordinate transformations}
\label{gauge}
\setcounter{equation}{0}

In this section, we show that there is a
choice of a local coordinate system $x^\alpha$
in which $u=x^0$ and the metric has the form
\bea
\lefteqn{
\metric = 2\zero{g}_{01}dx^0dx^1 +\zero{g}_{ab}dx^a dx^b
} \nonumber \\ 
&&
+\eps\left\{2\one{g}_{1a} dx^1 dx^a+\one{g}_{ab} dx^a
dx^b\right\}+O(\eps^2).
\label{rosen}
\eea
Here and below, indices $a,b,c,\dots$
take on the values $2,3$, while 
indices $i,j,k,\dots$ take on the values 
$1,2,3$.

The corresponding expansion of the contravariant form of
the metric tensor is
\bea
\lefteqn{
\metric^\sharp
=2\zero{g}{}\!^{01}\partial_0\partial_1
+\zero{g}{}\!^{ab}\partial_a\partial_b
} \nonumber
\\
&&
-\eps\left\{
2\zero{g}{}\!^{01}\zero{g}{}\!^{ab}\one{g}_{1b}\partial_0\partial_a
+\zero{g}{}\!^{ac}\zero{g}{}\!^{bd}
\one{g}_{cd}\partial_a\partial_b\right\}+O(\eps^2).
\label{contra}
\eea
For this metric, we have
\be
\zero{g}{}\!^{00} = 0, \qquad \one{g}{}\!^{00} = 0.\label{null3}
\ee
Thus, the phase $u=x^0$ satisfies \eq{null1}, and hence \eq{null}.

The most general coordinate transformation which is compatible with 
the expansion \eq{exp1} has the form
\bea
\frac{x^0}{\eps}
& \rightarrow &
\one{\Psi}{}\!^0\left(\frac{x^0}{\eps},x\right)
+\eps \two{\Psi}{}\!^0\left(\frac{x^0}{\eps},x\right)+O(\eps^2),
\label{t0}
\\
x^i 
& \rightarrow &
\zero{\Psi}{}\!^i(x)
+\eps \one{\Psi}{}\!^i\left(\frac{x^0}{\eps},x\right)
+\eps^2 \two{\Psi}{}\!^i\left(\frac{x^0}{\eps},x\right)+O(\eps^3).
\label{ti}
\eea
We suppose that the phase is given by $u=x^0$
in both the old and the new coordinates. 
Thus, the change of coordinates \eq{t0} implies a 
change in the phase of the form \eq{newu}.

First we simplify the leading order
metric components by means of a transformation
of the form
\be
x^0 \rightarrow x^0,
\qquad
x^i \rightarrow 
x^i+\eps\one{\Psi}{}\!^i\left(\frac{x^0}{\eps},x\right).
\label{trans1}
\ee
Expansion of the transformation
law for the change in covariant tensor components implies
that the leading order
metric components transform under \eq{trans1} according to
\bea
\zero{g}_{00}
& \rightarrow &
\zero{g}_{00}
+2 \one{\Psi}{}\!^k_{,\theta} \zero{g}_{0k}
+\one{\Psi}{}\!^k_{,\theta}\one{\Psi}{}\!^l_{,\theta} \zero{g}_{kl},
\label{t000}
\\
\zero{g}_{0i}
& \rightarrow &
\zero{g}_{0i}
+\one{\Psi}{}\!^k_{,\theta} \zero{g}_{ki},
\label{t00i}
\\
\zero{g}_{ij}
& \rightarrow &
\zero{g}_{ij}.
\label{t0ij}
\eea
If the matrix $\zero{g}_{ij}$ is non-singular,
then \eq{t00i} implies that we can transform $\zero{g}_{0i}$ to zero.
This contradicts the requirement that $x^0$ is null
(cf. \cite{LL}, Section 109). Hence, we must have
\be
\det \zero{g}_{ij}=0.
\label{det1}
\ee
By an appropriate renumbering of the $i$-coordinates, we
can suppose without loss of generality that
\be
\det \zero{g}_{ab}\ne 0.
\label{det2}
\ee

>From \eq{t000} and \eq{t00i}, we can then choose the 
transformation \eq{trans1} so that
\be
\zero{g}_{00}=\zero{g}_{02}=\zero{g}_{03}=0.
\label{goo}
\ee
Solving equation \eq{det1} for $\zero{g}_{11}$, we get
\be
\zero{g}_{11}=
\zero{g}{}\!^{ab}\zero{g}_{1a}\zero{g}_{1b},
\label{0g11}
\ee
where $\zero{g}{}\!^{ab}$ is the inverse
of $\zero{g}_{ab}$. We define
\be
g^a=\zero{g}{}\!^{ab}\zero{g}_{1b}.
\label{ga}
\ee
>From \eq{0g11}--\eq{ga}, it follows that
\be
\zero{g}_{11}=\zero{g}_{cd}g^cg^d,
\qquad
\zero{g}_{1a}=\zero{g}_{ac}g^c.
\label{0g11g1a}
\ee
Using \eq{goo}--\eq{0g11g1a} in \eq{metriccomp},
we find that in the transformed coordinate
system, the metric has the form
\be
\metric=2\zero{g}_{01}dx^0dx^1
+\zero{g}_{ab}\left(dx^a+g^a dx^1\right)
\left(dx^b+g^b dx^1\right)+O(\eps).
\label{line1}
\ee

>From \eq{pert1}, the metric \eq{line1} must satisfy
the condition
\be
\mtwo{R}_{ab}=0.
\label{rab=0}
\ee
Using \eq{line1} in \eq{connricciexp}, we find that
\be
\mtwo{R}_{ab}
=-\frac{1}{2}(\zero{g}{}\!^{01})^{2}\zero{g}_{ac}g^c_{,\theta}
\zero{g}_{bd}g^d_{,\theta}.
\label{rrrab}
\ee
Equations  \eq{rab=0}--\eq{rrrab} imply that
\[
g^a_{,\theta}=0,
\]
so $g^a$ is independent of $\theta$.
This fact allows us to remove $g^a$ by
a transformation
\be
x^a\rightarrow \Psi^a(x^1,x^c).
\label{trans-ga}
\ee
The form of the metric \eq{line1} is unchanged by
\eq{trans-ga}, and
\bea
\zero{g}_{ab}
&\rightarrow&
\Psi^c_{,a}\Psi^d_{,b}\zero{g}_{cd},
\nonumber
\\
g^a
&\rightarrow&
(A^{-1})^a_c (g^c+\Psi^c_{,1}),
\label{ga-trans}
\eea
where $(A^a_c)=(\Psi^a_{,c})$.
>From \eq{ga-trans}, we can  set $g^a=0$.
The metric \eq{line1} then reduces to
\be
\metric=2\zero{g}_{01}dx^0dx^1
+\zero{g}_{ab}dx^a dx^b+O(\eps).
\label{line1bis}
\ee

Next, we simplify the form of $\one{\metric}$.
We consider the transformation of coordinates
\be
x^0
\rightarrow
\eps \one{\Psi}{}\!^0\left(\frac{x^0}{\eps},x\right)
+\eps^2 \two{\Psi}{}\!^0\left(\frac{x^0}{\eps},x\right),
\quad
x^i
\rightarrow 
x^i+\eps^2 \two{\Psi}{}\!^i\left(\frac{x^0}{\eps},x\right).
\label{trans2}
\ee
Under the action of \eq{trans2}, the form \eq{line1bis}
of the metric is unchanged at order zero and
\[
\zero{g}_{01}\rightarrow 
\one{\Psi}{}\!^0_{,\theta}\zero{g}_{01}.
\]
At order one, the components transform according to
\beaa
\one{g}_{00}
& \rightarrow &
\one{\Psi}{}\!^0_{,\theta}
(\one{\Psi}{}\!^0_{,\theta}\one{g}_{00}
+2 \two{\Psi}{}\!^1_{,\theta}\zero{g}_{01}),
\\
\one{g}_{01}
& \rightarrow &
\one{\Psi}{}\!^0_{,\theta}\one{g}_{01}
+(\one{\Psi}{}\!^0_{,0}+\two{\Psi}{}\!^0_{,\theta})\zero{g}_{01},
\\
\one{g}_{0a}
& \rightarrow &
\one{\Psi}{}\!^0_{,\theta}\one{g}_{0a}
+\two{\Psi}{}\!^b_{,\theta}
\zero{g}_{ab},
\\
\one{g}_{11}
& \rightarrow &
\one{g}_{11}+2\one{\Psi}{}\!^0_{,1}\zero{g}_{01},
\\
\one{g}_{1a}
& \rightarrow &
\one{g}_{1a}+\one{\Psi}{}\!^0_{,a}\zero{g}_{01},
\\
\one{g}_{ab}
& \rightarrow &
\one{g}_{ab}.
\eeaa
These transformations can be used to make
\be
\one{g}_{11}=\one{g}_{0\alpha}=0.
\ee
The resulting metric then has the form given in \eq{rosen}.

Use of \eq{rosen} and \eq{contra} in \eq{connricciexp}
implies that the nonzero connection coefficients
at the orders $\eps^{-1}$ and $\eps^0$ are
\beaa
&&
\mone{\Gamma}{}\!^0{}_{00}
=\zero{g}{}\!^{01}\zero{g}_{01,\theta},
\quad
\mone{\Gamma}{}\!^1{}_{ab}
=-\frac{1}{2}\zero{g}{}\!^{01}\zero{g}_{ab,\theta},
\quad
\mone{\Gamma}{}\!^a{}_{0b}
=\frac{1}{2}\zero{g}{}\!^{ac}\zero{g}_{bc,\theta},
\\
&&
\zero{\Gamma}{}\!^0{}_{00}
=\zero{g}{}\!^{01}\zero{g}_{01,0},
\quad
\zero{\Gamma}{}\!^0{}_{0a}
=\frac{1}{2}\zero{g}{}\!^{01}
(\zero{g}_{01,a}+\one{g}_{1a,\theta})
-\frac{1}{2}\one{g}{}\!^{0b}\zero{g}_{ab,\theta},
\\
&&
\zero{\Gamma}{}\!^0{}_{ab}
=-\frac{1}{2}\zero{g}{}\!^{01}\zero{g}_{ab,1},
\quad
\zero{\Gamma}{}\!^1{}_{11}
=\zero{g}{}\!^{01}\zero{g}_{01,1},
\quad
\zero{\Gamma}{}\!^1{}_{1a}
=\frac{1}{2}\zero{g}{}\!^{01}
(\zero{g}_{01,a}-\one{g}_{1a,\theta}),
\\
&&
\zero{\Gamma}{}\!^1{}_{ab}
=-\frac{1}{2}\zero{g}{}\!^{01}
(\zero{g}_{ab,0}+\one{g}_{ab,\theta}),
\quad
\zero{\Gamma}{}\!^a{}_{01}
=-\frac{1}{2}\zero{g}{}\!^{ac}
(\zero{g}_{01,c}-\one{g}_{1c,\theta}),
\\
&&
\zero{\Gamma}{}\!^a{}_{0b}
=\frac{1}{2}\zero{g}{}\!^{ac}
(\zero{g}_{bc,0}+\one{g}_{bc,\theta})
-\frac{1}{2}\one{g}{}\!^{ac}\zero{g}_{bc,\theta},
\\
&&
\zero{\Gamma}{}\!^a{}_{1b}
=\frac{1}{2}\zero{g}{}\!^{ac}\zero{g}_{bc,1},
\quad
\zero{\Gamma}{}\!^a{}_{bc}
=\frac{1}{2}\zero{g}{}\!^{ad}
(\zero{g}_{bd,c}+\zero{g}_{cd,b}-\zero{g}_{bc,d})
+\frac{1}{2}\one{g}{}\!^{0a}\zero{g}_{bc,\theta}.
\eeaa
The nonzero components of the Ricci curvature at
the orders $\eps^{-2}$ and $\eps^{-1}$ are
\bea
&&
\mtwo{R}_{00}
=-\frac{1}{2}(\zero{g}{}\!^{ab}
\zero{g}_{ab,\theta})_{,\theta}
-\frac{1}{4}\zero{g}{}\!^{ac} \zero{g}_{bc,\theta}
\zero{g}{}\!^{bd} \zero{g}_{ad,\theta}
\nonumber
\\
&&
\qquad
+\frac{1}{2}\zero{g}{}\!^{01}\zero{g}_{01,\theta}
\zero{g}{}\!^{ab} \zero{g}_{ab,\theta},
\label{eq1}
\\
&&
\mone{R}_{01}
=-(\zero{g}{}\!^{01}\zero{g}_{01,1})_\theta
-\frac12(\zero{g}{}\!^{ab}\zero{g}_{ab,1})_{\theta}
-\frac14\zero{g}{}\!^{ac}\zero{g}_{bc,\theta}
\zero{g}{}\!^{bd}\zero{g}_{ad,1},
\label{R01-0}
\\
&&
\mone{R}_{ab}
=-\zero{g}{}\!^{01}\left(\zero{g}_{ab,1\theta}
-\frac{1}{2}\zero{g}{}\!^{cd}(\zero{g}_{ac,\theta}\zero{g}_{bd,1}
+\zero{g}_{ac,1}\zero{g}_{bd,\theta})\right.
\nonumber
\\
&&
\qquad
\left.
+\frac{1}{4}\zero{g}{}\!^{cd}(\zero{g}_{cd,1}\zero{g}_{ab,\theta}
+\zero{g}_{cd,\theta}\zero{g}_{ab,1})\right).
\label{Rab-0}
\\
&&
\mone{R}_{00}
=-\frac{1}{2}\one{g}{}\!^{a}_{a,\theta\theta}
-\frac{1}{2}\zero{g}{}\!^{ac}\zero{g}_{bc,\theta}\one{g}{}\!^{b}_{a,\theta}
+\frac{1}{2}\zero{g}{}\!^{01}\zero{g}_{01,\theta}\one{g}{}\!^{a}_{a,\theta}
-(\zero{g}{}\!^{ab}\zero{g}_{ab,\theta})_0
\nonumber
\\
&&
\qquad
-\frac{1}{2}\zero{g}{}\!^{ac}\zero{g}_{bc,\theta}
\zero{g}{}\!^{bd}\zero{g}_{ad,0}
+\frac{1}{2}\zero{g}{}\!^{01}\zero{g}{}\!^{ab}
(\zero{g}_{01,\theta}\zero{g}_{ab,0}
+\zero{g}_{01,0}\zero{g}_{ab,\theta})
\label{R00-0}
\\
&&
\mone{R}_{0a}
=\frac12(\zero{g}_{ab}\zero{g}{}\!^{01}
(\zero{g}_{01}\one{g}{}\!^{0b})_{\theta})_{\theta}
+\frac14\zero{g}{}\!^{cd}\zero{g}_{cd,\theta}
\zero{g}_{ab}\zero{g}{}\!^{01} (\zero{g}_{01}\one{g}{}\!^{0b})_{\theta}
\nonumber
\\
&&
\qquad
-\frac12(\zero{g}{}\!^{01}\zero{g}_{01,a}
+\zero{g}{}\!^{cd}\zero{g}_{cd,a})_\theta
+\frac12(\zero{g}{}\!^{bc}\zero{g}_{ab,\theta})_{c}
+\frac14\zero{g}{}\!^{bc}\zero{g}_{ab,\theta}\zero{g}{}\!^{de}\zero{g}_{de,c}
\nonumber
\\
&&
\qquad
+\frac14\zero{g}{}\!^{01}\zero{g}_{01,a}
\zero{g}{}\!^{cd}\zero{g}_{cd,\theta}
-\frac{1}{4}\zero{g}{}\!^{bd}\zero{g}_{cd,\theta}
\zero{g}{}\!^{ce}\zero{g}_{be,a},
\label{R0a-0}
\eea

\section{The asymptotic expansion}
\label{equations}
\setcounter{equation}{0}

We choose coordinates
\be
(x^0,x^1,x^2,x^3) = (u,v,y,z)
\label{uvyz}
\ee
in which the metric has the form \eq{rosen}.
We introduce functions $M$, $U$, $V$, $W$ 
of $(\theta,v,y,z)$ such that
\bea
\zero{g}_{01}&=&-e^{-M},
\nonumber\\
\left(\zero{g}_{ab}\right)
&=&\left(\begin{array}{cc}
e^{-U+V}\cosh W & -e^{-U}\sinh W \\
-e^{-U}\sinh W & e^{-U-V}\cosh W
\end{array}\right).
\label{lineernst}
\eea
It follows from \eq{rosen}, \eq{uvyz}, and \eq{lineernst}
that the leading order metric has the form of the colliding
plane wave metric,
\bea
&&\zero\metric = -2e^{-M}du\,dv+e^{-U}\bigl(e^{V}\cosh W dy^2
-2\sinh W dy\,dz
\nonumber
\\
&&\qquad\qquad\qquad\qquad\qquad\qquad\quad
+e^{-V}\cosh W dz^2\bigr).
\label{Brp}
\eea

>From \eq{eq1}, the only component of the
leading order perturbation equation \eq{pert1}
which is not identically satisfied is
\be
\mtwo{R}_{00}=0.
\label{r00}
\ee
Using \eq{eq1} and \eq{lineernst} in \eq{r00},
we obtain the $\theta$-constraint equation,
\be
U_{\theta\theta} = \frac12\left(U_{\theta}^2 + V_{\theta}^2\cosh^2 W
+ W_{\theta}^2\right)-U_{\theta}M_{\theta}.
\label{np0}
\ee

>From \eq{R01-0}--\eq{R0a-0}, the only components of
the first order perturbation equation \eq{pert2}
which are not identically satisfied are
\bea
&&
\mone{R}_{01}=0,\quad
\mone{R}_{ab}=0,
\label{rab}
\\
&&
\mone{R}_{00}=0, \quad \mone{R}_{0a}=0.
\label{r000b}
\eea
Using \eq{R01-0}--\eq{Rab-0} and \eq{lineernst} in \eq{rab},
we get the evolution equations in the colliding plane wave
equations,
\bea
U_{\theta v} &=& U_{\theta}U_{v},
\label{np1}
\\
V_{\theta v}
&=& \frac{1}{2}\left(U_{\theta}V_{v}+U_{v}V_{\theta}\right)
-\left( V_{\theta}W_{v}+V_{v}W_{\theta}\right)\tanh W,
\label{np2}
\\
W_{\theta v}
&=& \frac{1}{2}\left(U_{\theta}W_{v}+U_{v}W_{\theta}\right)
+V_{\theta}V_{v}\sinh W\cosh W.
\label{np3}
\\
M_{\theta v}&=&\frac12 \left(-U_{\theta} U_{v}
+ V_{\theta} V_{v}\cosh^2 W
+W_{\theta}W_{v}\right).
\label{np4}
\eea
>From \eq{contra} and \eq{R00-0}--\eq{R0a-0},
we find that \eq{r000b} is satisfied by a suitable
choice of the first order metric components $\one{g}_{ab}$,
$\one{g}_{1a}$.

\section{Variational principle}
\label{variational}
\setcounter{equation}{0}

The variational principle for the vacuum Einstein field equations is
\bea
&&\delta S = 0,\qquad S = \int L\, d^4x,\nonumber\\
&&L = R\sqrt{-\det g} ,\label{lagrangian}
\eea
where $R$ is the scalar curvature,
\[
R=g^{\alpha\beta}R_{\alpha\beta}.
\]

Using \eq{exp1}, \eq{exp2}, and \eq{Rexp} to expand the 
scalar curvature, we obtain
\bea
&&R=\frac{1}{\eps^2}\mtwo{R}
+\frac{1}{\eps}\mone{R}+O(1),
\nonumber\\
&&\mtwo{R}
=
\zero{g}{}\!^{\alpha\beta}\mtwo{R}_{\alpha\beta},
\label{R}
\\
&&\mone{R}
=
\zero{g}{}\!^{\alpha\beta}\mone{R}_{\alpha\beta}
-\one{g}{}\!^{\alpha\beta}\mtwo{R}_{\alpha\beta}.\nonumber
\eea
For a metric of the form \eq{line1bis},
we find that
\bea
&&\mtwo{R} = 0,\nonumber\\
&&\mone{R}
=\zero{g}{}\!^{ab}\mone{R}_{ab}
+2\zero{g}{}\!^{01}\mone{R}_{01} - \one{g}\!^{00}\mtwo{R}_{00}.
\label{scalarexp}
\eea
The only order one metric component which appears in \eq{scalarexp} is
\[
\lambda = -\one{g}\!^{00}.
\]
In the derivation of the asymptotic equations, we used a coordinate
system in which $\lambda=0$ --- see \eq{null3}.
In the variational principle, $\lambda$ acts as a Lagrange multiplier
for the constraint equation, so we do not set it
to zero until after we take variations.

We use \eq{scalarexp} in \eq{lagrangian}, expand the result
with respect to $\eps$, and write the expanded Lagrangian
in terms of $\lambda$ and the
functions $M$, $U$, $V$, $W$, defined in \eq{lineernst}.
This gives
\[
L = \frac{1}{\eps} \mone{L} + O(1),
\]
with
\bea
&&\mone{L}
= 
\left\{-2M_{\theta v}-4U_{\theta v}+3U_{\theta}U_v
+V_\theta V_v \cosh^2W+W_\theta W_v
\right\}e^{-U}\nonumber\\
&&\qquad + \lambda\left\{U_{\theta\theta}
- \frac12\left(U_{\theta}^2 + V_{\theta}^2\cosh^2 W
+ W_{\theta}^2\right)+U_{\theta}M_{\theta}\right\}e^{-M-U}
\label{I-1}
\eea
We make the change of variables in the integration
\[
d^4x = du\, dv\, dy\, dz = \eps d\theta\, dv \, dy \, dz,
\]
and neglect the integration with respect to the parametric variables
$(y,z)$. The leading order asymptotic variational
principle then becomes
\[
\delta \zero{S} = 0,\qquad \zero{S} = \int \mone{L}\,d\theta dv.
\]
Variations of $\zero{S}$ with respect
to the first order metric component $\lambda$
give the constraint \eq{np0}.
Variations with respect to $M$, $U$, $V$, $W$
give the evolution equations \eq{np1}--\eq{np4},
after we set $\lambda=0$. It is permissible to set $\lambda = 0$
because the constraint is a gauge-type
constraint which is preserved by the evolution equations.

\section{Boundary conditions}
\label{jump}
\setcounter{equation}{0}

In this section, we discuss the derivation of boundary conditions
for the asymptotic equations. For simplicity, we consider a ``sandwich'' wave
located near the null surface $u = 0$ which varies rapidly in a thin strip
\[
\theta_- \le \frac{u}{\eps} \le \theta_+.
\]
We denote the slowly varying metrics on either side
of the wave by
\be
\metric = \left\{ \begin{array}{ll}
	\metric_+ &\mbox{in $u > 0$} \\
	\metric_- &\mbox{in $u < 0$}
	\end{array}\right. .\label{metricpm}
\ee

We consider a coordinate patch around a point on
the surface $u=0$ with local coordinates $(u,v,y,z)$ chosen
as in the derivation of the asymptotic solution.
In order for the metric outside the wave to join continuously
with the solution inside, we must have
\bea
&&\lefteqn{
\metric_\pm \to  -2e^{-M_\pm}du\,dv
+e^{-U_\pm}(e^{V_\pm}\cosh W_\pm dy^2
}
\nonumber
\\
&&\quad
-2\sinh W_\pm dy\,dz
+e^{-V_\pm}\cosh W_\pm dz^2),\label{Brpm}
\eea
as $u \to 0^{\pm}$, where $M_\pm$, $U_\pm$, $V_\pm$, $W_\pm$
are functions of $(v,y,z)$. From \eq{Brp}, \eq{Brpm}, and the
continuity of the metric, it follows that the solution
of  \eq{np1}--\eq{np4} must satisfy the characteristic boundary conditions,
\be
M = M_\pm, \quad U = U_\pm, \quad V = V_\pm, \quad 
W = W_\pm, \quad \mbox{when $\theta = \theta_\pm$}.
\label{charbc}
\ee
This data need not satisfy the constraint \eq{npc6}.

The asymptotic equations must be supplemented by a condition
which specifies the profile of the wave. For example,
we can impose a characteristic initial condition
\be
M = M_0, \quad U = U_0, \quad V = V_0, \quad 
W = W_0, \quad \mbox{when $v=0$},
\label{charic}
\ee
where
$M_0$, $U_0$, $V_0$, $W_0$ are functions of $(\theta,y,z)$
which satisfy the constraint \eq{np0}.
The characteristic initial data must also be compatible with
the characteristic boundary data, meaning that
\[
M_0(\theta_\pm,y,z) = M_\pm(0,y,z),
\]
together with analogous conditions for the other variables.

Equations \eq{np1}--\eq{np4}, the characteristic initial
condition \eq{charic} on $v=0$, and the characteristic boundary
condition \eq{charbc} on $\theta=\theta_-$ form a well-posed problem.
Provided that the solution inside the wave is free of singularities,
this problem has a unique solution. In particular, the solution
at $\theta=\theta_+$ is uniquely determined.
Thus, in principle, the asymptotic equations, together with the
characteristic initial data \eq{charic}, determine  a set of
jump relations which connect the minus and plus metrics ahead of
and behind the wave, respectively.
If the metric ahead of the wave is known, then the jump
conditions provide characteristic boundary conditions 
on $u=0$ for the space-time behind the wave. Together with a
characteristic initial condition on $v=0$ and $u > 0$, for example,
this gives a characteristic initial value
problem \cite{Re} for the full field equations.
This problem determines the slowly varying
metric behind the wave (at least locally).

For instance, in the case of a plane polarized wave, the solution
of \eq{apc1} for $U$ is \cite{Gr}
\be
U(\theta,v) = -\log\left[f(\theta) + g(v)\right].
\label{solU}
\ee
Here $f$ and $g$ are functions of integration, and
we do not explicitly show the possible parametric dependence
of the functions on $(y,z)$. From \eq{charbc}, \eq{charic},
and \eq{solU} we have
\[
f(\theta) + g(0) = e^{-U_0(\theta)},\qquad
f(\theta_-) + g(v) = e^{-U_-(v)}.
\]
The solution is nonsingular provided that
$f(\theta)+g(v) > 0$.

It follows from \eq{solU} that the jump in $U$ satisfies
\[
e^{-U_+(v)} - e^{-U_-(v)} = e^{-U_0(\theta_+)} - e^{-U_0(\theta_-)}.
\]
Use of \eq{solU} in \eq{apc2} gives a linear wave equation for $V$,
\[
\left(f+g\right)V_{\theta v} = \frac12\left(g_vV_\theta+ f_\theta V_v\right).
\]
Solution of this equation with the characteristic initial data
$V = V_0$ on $v=0$ and the characteristic boundary data
$V = V_-$ on $\theta = \theta_-$ determines, in principle,
the solution $V = V_+$ on $\theta=\theta_+$.
Finally, when $W=0$, we define the $v$-constraint function $G$ by
\be
G = U_{vv} - \frac{1}{2}\left(U_{v}^2+V_{v}^2\right)+U_{v}M_{v}.
\label{defG}
\ee
It follows from \eq{defG} and \eq{apc1}--\eq{apc3} that
\[
G_\theta = U_\theta G.
\]
Integration of this equation with respect to $\theta$ implies that
\[
\log G_+(v) - \log G_-(v) = U_+(v) - U_-(v).
\]
This equation provides a jump condition for $M$.

One difficulty which arises in the formulation of boundary conditions
ahead of the wave is that the metric $\metric_-$ may not be given
in a coordinate system which is compatible with the coordinate
system used in the derivation of the asymptotic equations. It is then
necessary to construct compatible coordinates $(u,v,y,z)$.
The $u$-coordinate is the phase, so it is a null coordinate
of the metric which can be found by solving an eikonal
equation, subject to appropriate initial conditions. 
The $v$-coordinate is a null coordinate which is
orthogonal to $u$, while the $y$ and $z$ coordinates parametrize
the null geodesics on the surface $u=v=0$.

If the gravitational wavefront $u=0$ forms a caustic, then
the solution of the eikonal equation becomes multi-valued.
When this happens, the local plane-wave approximation breaks
down, and the asymptotic solution is not valid. However,
the focusing at a caustic of the congruence of null geodesics
associated with the phase does not necessarily imply the formation
of a space-time singularity.

\section{Examples}
\label{examples}
\setcounter{equation}{0}

In this Section, we derive boundary conditions
for the asymptotic equations 
which correspond to the propagation of a non-planar
gravitational wave into Minkowski space-time, the exterior
Schwarzschild space-time, and Robertson-Walker space-time.
In each example, we consider the case of spherical
waves, where the boundary data can be explicitly computed.
In this paper, we do not attempt to explore the physical
consequences of the asymptotic equations in any detail.
Our aim here is simply to illustrate how to apply the
asymptotic equations to specific physical problems.

\subsection{Nonplanar wave propagation into Minkowski space-time}

We suppose that the space-time ahead of the wave is flat.
In inertial coordinates $(t,\vec x)$,
with $t = x^0$ and $\vec x = (x^1,x^2,x^3)$, the metric is
\[
\metric_- = - dt^2 + d\vec x^2.
\]
We consider a wave with phase
\[
u = \frac{t - w(\vec x)}{\sqrt{2}}.
\]
The phase $u$ is a null function of $\metric_-$ if
\[
|\nabla w|^2 = 1,
\]
where $\nabla$ is the gradient with respect to $\vec x$.
We define
\[ 
v =  \frac{t + w(\vec x)}{\sqrt{2}},
\]
and choose coordinates $y(\vec x)$, $z(\vec x)$ such that
$\nabla w$, $\nabla y$, $\nabla z$ are orthogonal.
In the $(u,v,y,z)$ coordinates, we have
\be
\metric_- = - 2 du\,dv  + \frac{1}{|\nabla y|^2}dy^2 +
\frac{1}{|\nabla z|^2}dz^2.
\label{flatuv}
\ee
A comparison of \eq{flatuv} with \eq{Brpm} shows that
the minus boundary data is given by
\[
M_- = 0,
\quad e^{-U_-} = \frac{1}{\left.|\nabla y| |\nabla z|\right|_{u=0}},
\quad e^{-V_-} = \left.\frac{|\nabla y|}{|\nabla z|}\right|_{u=0},
\quad W_- = 0.
\]

For example, in the case of an outgoing spherical wave,
suitable coordinates are
\be
u = \frac{t - r}{\sqrt{2}},\quad
v =  \frac{t + r}{\sqrt{2}},\quad
y = \tea,\quad
z = \vphi,\label{sphco}
\ee
where $(r,\tea,\vphi)$ are spherical polar coordinates and $t > 0$.
In $(u,v,y,z)$ coordinates, the flat space-time metric is
\[
\metric_- = -2 du dv + \frac{1}{2} (u-v)^2
\left(dy^2 + \sin^2y\, dz^2\right).
\]
Evaluation of this metric
at $u=0$ and a comparison with \eq{Brpm}
gives the minus boundary data
\[
M_- = 0,
\quad e^{-U_-} = \frac12 v^2 \sin y,
\quad e^{-V_-} = \sin y,
\quad W_- = 0,
\]
where $v>0$.
In this case, $M_-$, $V_-$, and $W_-$ are
independent of $v$, while $U=U_-$
satisfies the equation
\[
U_{vv} = \frac{1}{2} U_{v}^2.
\]
Thus, the boundary data satisfies the $v$-constraint equation \eq{npc6}.
The solution is therefore identical to an exact solution for the collision
of outgoing and incoming spherical waves, with the additional
possibility of a slow parametric dependence on the polar angles $(\tea,\vphi)$.
Some exact solutions for spherical wave propagation into flat space-time
are constructed in \cite{Al}.

For an incoming spherical wave, we use
\[
u = \frac{t + r}{\sqrt{2}},\quad
v =  \frac{t - r}{\sqrt{2}},
\]
where $t < 0$. This leads to the same boundary data as in the case
of an outgoing spherical wave, but with $v < 0$ instead of $v >0$.

\subsection{Gravitational waves incident on a black hole}

The exterior Schwarzschild metric is
\be
\metric_- = -a dt^2 + \frac{1}{a} dr^2 +
r^2 \left(d\tea^2 + \sin^2\tea d\vphi^2\right),
\label{scmetric}
\ee
where $r > 2m$ and
\[
a(r) = 1 - \frac{2m}{r}.
\]
The contravariant metric tensor is
\[
\metric_-^\sharp = -\frac{1}{a}\partial_t^2 + {a} \partial_r^2 +
\frac{1}{r^2} \left(\partial_\tea^2 +
\frac{1}{\sin^2\tea} \partial_\vphi^2\right).
\]
For simplicity, we consider an axially symmetric phase of the form
\[
u = \frac{t - w(r,\tea)}{\sqrt{2}}.
\]
The function $u$ is null if $w$ satisfies the eikonal equation
\[
a w_r^2 + \frac{1}{r^2} w_\tea^2 = \frac{1}{a}.
\]
We define the orthogonal null coordinate
\[
v = \frac{t + w(r,\tea)}{\sqrt{2}},
\]
and choose a coordinate $y(r,\tea)$ whose gradient is orthogonal
to the gradient of $w(r,\tea)$. In that case, we have
\[
y_r = -\frac{h w_\tea}{r^2},\qquad y_\tea = h a  w_r,
\]
where $h(r,\tea)$ is a suitable integrating factor. We take
$z = \vphi$. In $(u,v,y,z)$ coordinates, the Schwarzschild metric
\eq{scmetric} is given by
\be
\metric_- = - 2a du\,dv + \frac{r^2}{h^2} dy^2 +
r^2 \sin^2\tea dz^2.
\label{scuv1}
\ee
Inversion of the change of coordinates
$(t,r,\tea,\vphi) \mapsto (u,v,y,z)$
implies that $r=r_-(v,y)$ and $\tea = \tea_-(v,y)$ on $u=0$
for suitable functions $r_-$ and $\tea_-$. A comparison
of \eq{scuv1} with \eq{Brpm} implies that
the boundary data is given by
\[
e^{-M_-} = a_-,\quad e^{-U_-} = \frac{r_-^2\sin\tea_-}{h_-},
\quad e^{-V_-} = h_- \sin\tea_-, \quad W_- = 0,
\]
where $a_-=a(r_-)$ and $h_- = h(r_-, \tea_-)$.

In the case of an incoming spherical wave incident on
the black hole, suitable coordinates are
\be
u = \frac{t + A(r)}{\sqrt{2}},\quad
v =  \frac{t - A(r)}{\sqrt{2}},\quad
y = \tea,\quad
z = \vphi,\label{sphbh}
\ee
where
\[
A_r = \frac{1}{a}.
\]
Integration of this equation implies that
\[
A(r) = r + \log(r-2m).
\] 
In $(u,v,y,z)$ coordinates, the exterior Schwarzschild metric is
\be
\metric_- = -2a \, du dv +
r^2 \left(dy^2 + \sin^2 y\, dz^2\right).
\label{scuv}
\ee
>From \eq{sphbh}, we have $r = r_-(v)$ on $u=0$
where the function $r_-(v)$ is the solution of
\be
A(r_-) = -\frac{v}{\sqrt{2}}.\label{defrminus}
\ee
A comparison of \eq{scuv} with \eq{Brpm} implies
that the boundary data ahead of
the incoming spherical wave is given by
\be
e^{-M_-} = a_-,\quad e^{-U_-} = r_-^2\sin y,
\quad e^{-V_-} = \sin y, \quad W_- = 0.
\label{scbc}
\ee
Dropping the minus subscripts, we find that
the constraint function $G$ in \eq{defG} for the boundary
data \eq{scbc} is given by
\[
G = 2\left(\frac{a_vr_v}{ar} - \frac{r_{vv}}{r}\right).
\]
Differentation of \eq{defrminus} with respect to $v$ implies that
\[
r_v = -\frac{a}{\sqrt2},\qquad r_{vv} = -\frac{a_v}{\sqrt2}.
\]
Use of this equation in the expression for $G$ implies that $G = 0$.
Thus, the boundary data \eq{scbc} satisfies the $v$-constraint equation
\eq{npc6}.

Numerical solutions of the interaction of a spherical gravitational wave with
a black hole appear in \cite{Sm}.

\subsection{Gravitational waves in a Robertson-Walker space-time}

The Robertson-Walker metric is
\be
\metric_- = - dt^2 + \frac{1}{R^2} \left\{\frac{1}{1-kr^2}dr^2 +
r^2 \left(d\tea^2 + \sin^2\tea d\vphi^2\right)\right\},
\label{rwmetric}
\ee
where $R(t)$ is the scale factor, and $k = -1,0,1$. 

As in the Schwarzschild example, we consider an axially symmetric phase
for simplicity, given by
\[
u = \frac{I(t) - w(r,\tea)}{\sqrt{2}},
\]
where
\be
I_t = R,\qquad
\left(1-kr^2\right)w_r^2 + \frac{1}{r^2} w_\tea^2 = 1.
\label{It}
\ee
We define an orthogonal null coordinate $v$ by
\[
v = \frac{I(t) + w(r,\tea)}{\sqrt{2}}.
\]
We choose a coordinate $y(r,\tea)$ whose gradient is orthogonal
to the gradient of $w(r,\tea)$, so that
\[
y_r = -\frac{h w_\tea}{r^2\sqrt{1-kr^2}},\qquad y_\tea = h \sqrt{1-kr^2} w_r,
\]
where $h(r,\tea)$ is a suitable integrating factor, and take
$z = \vphi$. In $(u,v,y,z)$ coordinates, the Robertson-Walker metric
\eq{rwmetric} is given by
\be
\metric_- = - \frac{2}{R^2} du\,dv + \frac{r^2}{h^2R^2} dy^2 +
\frac{r^2}{R^2} \sin^2\tea dz^2.
\label{rwuv}
\ee
A comparison of \eq{rwuv} with \eq{Brpm}
implies that the boundary data is given by
\[
e^{-M_-} = \frac{1}{R_-^2},
\quad e^{-U_-} = \frac{r_-^2\sin\tea_-}{h_-R_-^2},
\quad e^{-V_-} = h_- \sin\tea_-, \quad W_- = 0,
\]
where $r=r_-(v,y)$, $\tea = \tea_-(v,y)$, $R=R_-(v,y)$,
and $h = h_-(v,y)$ on $u=0$.

For an outgoing spherical wave in a Robertson-Walker space-time,
suitable coordinates are
\[
u = \frac{I(t) - w(r)}{\sqrt{2}},\quad
v =  \frac{I(t) + w(r)}{\sqrt{2}},\quad
y = \tea,\quad
z = \vphi,
\]
where
\[
w_r = \frac{1}{\sqrt{1-kr^2}}.
\]
Integration of this equation implies that
\[
w(r) = \left\{
\begin{array}{ll}
\sin^{-1}r & \mbox{if $k=1$}, \\
r & \mbox{if $k=0$}, \\
\sinh^{-1}r & \mbox{if $k=-1$}.
\end{array}\right.
\] 
The corresponding boundary data is given by
\be
e^{-M_-} = \frac{1}{R^2_-},\quad e^{-U_-} =
\frac{r_-^2 \sin y}{R_-^2},
\quad e^{-V_-} = \sin y, \quad W_- = 0,\label{rwbd}
\ee
where $t_-(v)$ and $r_-(v)$ are
given by
\be
I(t_-) = \frac{v}{\sqrt 2},\qquad
r_- = \left\{
\begin{array}{ll}
\sin(v/\sqrt2) & \mbox{if $k=1$}, \\
v/\sqrt2 & \mbox{if $k=0$}, \\
\sinh(v/\sqrt2) & \mbox{if $k=-1$},
\end{array}\right.
\label{Iw}
\ee
and  $R_- = R(t_-)$.

Dropping the minus subscripts, we find that
the constraint function $G$ in \eq{defG}
for the boundary data \eq{rwbd} is given by
\be
G = 2\left(\frac{R_{vv}}{R} - \frac{r_{vv}}{r}\right).
\label{Grw}
\ee
>From \eq{It} and \eq{Iw}, we find that
\[
r_{vv}= -\frac12 k r,\quad R_{vv} = \frac{RR_{tt}-R_t^2}{2R^3}.
\]
Use of these expressions in \eq{Grw} gives
\[
G = \frac{RR_{tt}-R_t^2}{R^4} + k.
\]
Thus, in general, the boundary data \eq{rwbd}
does not satisfy the $v$-constraint
equation \eq{npc6}.

\bigskip\bigskip\noindent
{\bf Acknowledgments.} The work of J.K.H. was partially supported by the
NSF under grant number DMS-9704152.
The work of G.A. was partially supported by the CNR
1996/97 short-term fellowships program.

\end{document}